\documentclass[12pt]{article}
\usepackage{geometry}                
\usepackage[parfill]{parskip}    
\usepackage{graphicx}
\usepackage{amssymb}
\usepackage{amsthm}
\usepackage{amsmath}
\usepackage{amsfonts}
\usepackage{epstopdf}
\usepackage{url}
\theoremstyle{plain}
\newtheorem{thm}{Theorem}[section]
\newtheorem{lem}[thm]{Lemma}
\newtheorem{cor}[thm]{Corollary}

\newtheorem{conj}[thm]{Conjecture}

\theoremstyle{definition}
\newtheorem{defn}[thm]{Definition}
\newtheorem{q}[thm]{Question}

\begin{document}

\title{Zero Forcing Sets and Bipartite Circulants}

 \author{Seth A. Meyer\\
 Department of Mathematics\\
 University of Wisconsin\\
 Madison, WI 53706\\
 {\tt smeyer@math.wisc.edu}
 }

\maketitle

 \begin{abstract} 
In this paper we introduce a class of regular bipartite graphs whose biadjacency matrices are circulant matrices and we describe some of their properties.  Notably, we compute upper and lower bounds for the zero forcing number for such a graph based only on the parameters that describe its biadjacency matrix.  The main results of the paper characterize the bipartite circulant graphs that achieve equality in the lower bound.


\end{abstract}

  \begin{section}{Introduction}

  In this introduction, we give an overview of the paper and general motivation for each topic.  We refer the reader to specific sections later in the paper for formal definitions, results, and proofs.  In Section \ref{BC} we define a family of regular bipartite graphs which we call {\em bipartite circulants} that have properties similar to those possessed by circulant graphs and circulant digraphs, see e.g. \cite{davis}.  This allows us to consider biadjacency matrices (freeing us from the symmetric requirement for adjacency matrices) while still keeping several nice properties of circulant graphs.  It is then possible for us to define one of these graphs by its biadjacency matrix with just a few parameters and compute graph statistics based only on the parameters chosen.\\
  
  In Section \ref{ZF} we review the idea of the zero forcing number for a graph, $Z(G)$, as defined in \cite{AIM08}.  This is a graph statistic that can be computed surprisingly often for our bipartite circulants based on their construction.  We focus mainly on finding bounds for $Z(G)$ based only on the matrix parameters.  After obtaining a trivial lower bound on $Z(G)$ based the regularity of $G$ and the fact that it is bipartite, we then discuss the parameters for a bipartite circulant which achieve equality in the bound.\\
  
  In their full generality bipartite circulants can be quite complicated, so in Section \ref{cubic} we examine what can be said when we restrict ourselves to 3-regular bipartite circulants (we will refer to 3-regular graphs as cubic in the rest of the paper).  This is the first non-trivial case for bipartite circulants, and the upper and lower bounds from Section \ref{ZF} can be extended in this case.  Notably, we prove that the unique bipartite, cubic graph that obtains equality in the lower bound from the previous section is a bipartite circulant.\\
  
  Finally, we end in Section \ref{conj} with some open questions and potential directions of future research.
  \end{section}
  
  \begin{section}{Bipartite Circulant Graphs}\label{BC}
  
    Recall that an $m\times n$ matrix $A$ with entries equal to either zero or one (we will refer to these as (0,1) matrices in the rest of the paper) can be thought of as the biadjacency matrix of a bipartite graph on $m+n$ vertices.  In this correspondence, we have a vertex for each row and column in the matrix.  Two vertices are adjacent if and only if one is a row vertex, the other is a column vertex, and the corresponding entry in $A$ is a 1.  The $m\times n$ (0,1) matrices are thus in bijection with labeled bipartite graphs on $m+n$ vertices with bipartitions of size $m$ (the rows) and $n$ (the columns).  To visualize this, draw the row vertices in a vertical line on the left and the column vertices in a vertical line on the right, we refer to a specific vertex by its label and its type, e.g. the row vertex with label 3 with be denoted $3$L.  This section is concerned with graphs rather than matrices, and so we ignore the labels, which allows two different matrices to represent the same bipartite graph when they are attained from each other by independent row and column permutations.\\
    
    A square matrix is a circulant matrix if each row in the matrix is a copy of the row above it with all the entries shifted one position to the right, with the last entry wrapping to the first column.  If $A$ is an $n\times n$ (0,1) circulant matrix we can consider the bipartite graph generated by $A$.  This leads to our first definition.\\
    
    \begin{defn}A bipartite graph obtained from a square (0,1) circulant matrix $A$ is said to be a bipartite circulant graph, or more colloquially, a bipartite circulant.
    \end{defn}
    
    When it is convenient we will identify a bipartite graph $G$ with its biadjacency matrix $A$, we write this as $G\simeq A$.  We use the symbol $\simeq$ because a given bipartite circulant graph can usually be represented by several biadjacency matrices, often it can even be represented by more than one circulant matrix.  For example, the matrices\\
    
    \begin{center}
    $A=\left[\begin{array}{cccccc}
    1 & 1 & 0 & 1 & 0 & 0\\
    0 & 1 & 1 & 0 & 1 & 0\\
    0 & 0 & 1 & 1 & 0 & 1\\
    1 & 0 & 0 & 1 & 1 & 0\\
    0 & 1 & 0 & 0 & 1 & 1\\
    1 & 0 & 1 & 0 & 0 & 1\end{array}\right]$ and
    $A'=\left[\begin{array}{cccccc}
    0 & 1 & 1 & 0 & 1 & 0\\
    0 & 0 & 1 & 1 & 0 & 1\\
    1 & 0 & 0 & 1 & 1 & 0\\
    0 & 1 & 0 & 0 & 1 & 1\\
    1 & 0 & 1 & 0 & 0 & 1\\
    1 & 1 & 0 & 1 & 0 & 0\end{array}\right]$
    \end{center}
    \vspace{3 mm}
    both correspond to the following graph:\\

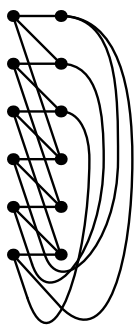
    \begin{center}
    \includegraphics{Graphics.2}
    \end{center}
    
    Any (0,1) circulant matrix can also be defined as a sum of distinct powers of $P$, the permutation matrix corresponding to the $n$-cycle $(1\ 2\ 3\ \cdots\ n)$.  For example, if $n=6$,\\
    
        \begin{center}
        $P= \left[\begin{array}{cccccc}
        0 & 1 & 0 & 0 & 0 & 0\\
        0 & 0 & 1 & 0 & 0 & 0\\
        0 & 0 & 0 & 1 & 0 & 0\\
        0 & 0 & 0 & 0 & 1 & 0\\
        0 & 0 & 0 & 0 & 0 & 1\\
        1 & 0 & 0 & 0 & 0 & 0\end{array}\right].$
        \end{center}
        
    We have $I_n+P^1+P^2+\cdots +P^{n-1}=J_n$, the $n\times n$ matrix of all 1s, so we only need to consider $0,1,\ldots,n-1$ as powers of $P$.  Note that $P^0=I_n$, the $n\times n$ identity matrix.  To standardize our notation we assume that if $A=P^{i_1}+P^{i_2}+\cdots+P^{i_k}$ we have $0\leq i_1<i_2<\cdots<i_k\leq n-1$.  With this convention the matrices above are written $A=I+P+P^3$ and $A'=P+P^2+P^4$.  Determining when two matrices with different powers of $P$ give bipartite circulants that are isomorphic will be a recurring problem for us as we move forward.\\
    
    We will need several basic properties of these graphs in order to continue our investigation.\\
    
    \begin{lem}\label{BCprops}
    Let $G$ be a connected bipartite circulant graph on $2n$ vertices with $G\simeq P^{i_1}+P^{i_2}+\cdots + P^{i_k}$.  Then the following are true:
      \begin{enumerate}
      \item $G$ is $k$-regular,
      \item We can refer to vertices by their position in a biadjacency matrix of $G$.  In this fashion $P^{i_j}$ corresponds to edges from $m$L to $(m+i_j)$R for all $m$.  In this characterization the neighbors of $m$L are $\{(m+i_1)R, \ldots, (m+i_k)R\}$ and the neighbors of $m$R are $\{(m-i_1)L,\ldots,(m-i_k)L\}$,
      \item Without loss of generality we may assume that $i_1=0$.
      \end{enumerate}
    \end{lem}
    \begin{proof}
    Since each row in the specified biadjacency matrix has exactly one 1 for each power of $P$ in the sum, each row vertex is adjacent to $k$ column vertices.  Similarly, each column has one 1 for each power of $P$ so every column vertex is adjacent to $k$ row vertices.  Thus, $G$ is $k$-regular.
    
    To see (2), note that for any given row, say row $j$, $P^{i_1}$ has a one in column $j+i_1$.  This gives us that $j$L is adjacent to $(j+i_1)$R for each $j$.  The sum of the powers of $P$ then specifies all of the edges in this fashion.
    
    To prove (3), we label the vertices of $G$ by their positions in the biadjacency matrix.  Then we relabel the column vertices by subtracting $i_1$ from each of their labels, from (2) we have that $j$L is adjacent to $j$R for each $j$.  This corresponds to cycling the columns of the matrix until we have brought a diagonal of 1s into the main diagonal position.  If the biadjacency matrix was a circulant, it's still a circulant after cycling the columns so we have the result.
    \end{proof}
    
    Most questions about graphs can be reduced to questions about their connected components, so we find it natural to describe when a bipartite circulant is connected.\\
    
    \begin{lem}\label{conn}
     Let $G$ be a bipartite circulant on $2n$ vertices with $G\simeq I+P^{i_2}+P^{i_3}+\cdots+P^{i_k}$.  Then $G$ is connected if and only if gcd$(i_2,\ldots,i_k,n)=1$.
    \end{lem}
    \begin{proof} Since we can assume without loss of generality that $I$ is a term in the biadjacency matrix from Lemma \ref{BCprops}, we see that vertex $i$L is adjacent to vertex $i$R in $G$. To prove the forward implication we note that if $G$ is connected then there is a path from 0L to 1L.  This path uses some right-to-left edges (corresponding to subtracting) and some left-to-right edges (corresponding to adding); say that we use $e_{j,1}$ right-to-left edges from $P^{i_j}$ and $e_{j,2}$ left-to-right edges from $P^{i_j}$.  Then we see that \[ \displaystyle \sum_{j=2}^{k}(e_{j,2}-e_{j,1})i_j \equiv 1 \bmod{n}.\]  This gives us a linear combination of $\{i_2,\ldots,i_k\}$ which is congruent to 1 modulo $n$.  Equivalently, we have a linear combination of $\{i_2,\ldots,i_k,n\}$ that equals 1, proving the implication.\\
    
    The other direction is similar.  If we have gcd$(i_2,\ldots,i_k,n)$=1 then we have a linear combination of the $i_j$ that determines a pattern of edges in the graph corresponding to traveling from $m$L to $(m+1)$L for any any $m$.  Combining this with the fact that $I$ is a term in the matrix gives us a path from $0$L to any other vertex, so $G$ is connected.
    \end{proof}
    
    It is easy to see that starting from 0L following alternating edges from $P^{j}$ and $I$ forms a cycle whose length is the order of $j$ in the additive group $\mathbb{Z}/n\mathbb{Z}$.  We know that the order of an element $j$ in the cyclic group of order $n$ is n/gcd$(j,n)=n/d$, and we can use this to draw a picture of such a graph in a more illuminating way.  Because the cosets of $\langle j \rangle$ all have the same size, we see that if $G'\subseteq G$ is the subgraph of $G$ with only edges from $I$ and $P^{i_j}$, $G'$ decomposes into $d$ cycles of length $\frac{2n}{d}$.  Since each cycle will correspond to a different coset of $\langle j \rangle$ in the group, we have one cycle with all the vertices congruent to zero mod $d$, another with all the vertices congruent to 1 mod $d$, and so on up to a cycle with those vertices congruent to $(d-1)$ mod $d$.  We can then observe what the other edges in $G$ are congruent to mod $d$, which determines how the cycles are connected by the edges not in $G'$.\\
    
    It is interesting to note that it is possible, even in the cubic case, to have both non-zero powers divide $n$, which leads to two different decompositions of $G$ into cycles.  $G\simeq I+P^3+P^5$ with $n=15$, pictured below illustrates this point.  The edges from $I$ go left-to-right in each of the three vertical 10-cycles, the edges of $P^5$ go left-to-right from one 10-cycle to the next, and the edges of $P^3$ connect one horizontal level with the adjacent levels.  We will discuss this case further in Section \ref{cubic}.\\
    
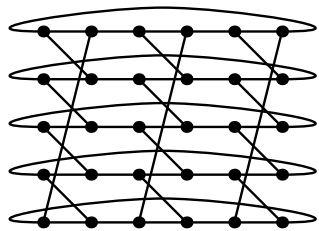
    \begin{center}
    \includegraphics{Graphics.1}
    \end{center}
    
    There are cycles of length 6 reading left to right across the top (0L, 0R, 10L, 10R, 5L, 5R) and cycles of length 10 reading down in the first column (0L, 3R, 3L, 6R, and so on).  $G$ is still connected because even though both gcd($5,15)=5>1$ and gcd$(3,15)=3>1$ we have gcd(3,5,15)=1, satisfying Lemma \ref{conn}.\\
    
    We now have two results that depend on the specific powers of $P$ present in the biadjacency matrix of $G$.  Recall that these are not necessarily uniquely defined: $I+P+P^3\simeq G\simeq P+P^2+P^4$ for $n=6$.  The next result gives conditions that allow us to generate circulant matrices that give isomorphic (unlabeled) bipartite graphs with different powers of $P$.\\
    
    \begin{thm}\label{BCiso}
    Let $G$ be a connected bipartite circulant on $2n$ vertices with $G\simeq P^{i_1}+\cdots +P^{i_k}$.  Then for every unit $u$ in the ring $\mathbb{Z}/n\mathbb{Z}$ and element $z\in \mathbb{Z}/n\mathbb{Z}$ we have $G\simeq P^{ui_1+z}+ \cdots + P^{ui_k+z}$.
    \end{thm}
    \begin{proof}
    We will prove the part of the theorem corresponding to adding $z$ first.  We know that if $A=P^{i_1}+\cdots + P^{i_k}$ is a biadjacency matrix of $G$ we can obtain other biadjacency matrices of $G$ by permuting the rows and/or columns of $A$.  This can be written as $G\simeq QAR$ where $Q$ and $R$ are permutation matrices.  Since $P$ is a permutation matrix, albeit a special one in our investigation, we have that
    \[G\simeq IAP^{z} = P^{i_1+z}+\cdots + P^{i_k+z},\]
    which is the desired result.  We note that this corresponds to adding $z$ to the labels of one of the sets of vertices in the bipartition of $G$.  This can be seen in the graph as moving the right side vertices down $z$ positions, which makes each edge from a left side vertex go an additional $z$ positions.\\
    
    To obtain multiplication, we relabel each vertex with $u$ times its original label.  This relabeling is bijective because $u$ is a unit, so this operation is well defined.  This changes an existing edge going from $m$L to $(m+i_l)$R into an edge going from $(mu)$L to $(m+i_l)u$R.  Since $(m+i_l)u=mu+i_l u$, edges from $P^{i_l}$ now correspond to adding $ui_l$, and we can write $G\simeq P^{ui_1}+\cdots P^{ui_k}$.  Combining these two procedures gives the result, since any sequence of multiplications by units and additions can be written as only one multiplication by a unit and one addition by distributing the multiplication.
    \end{proof}
    
    We now turn our attention to computing a graph statistic for bipartite circulants, the zero forcing number of a graph.
  \end{section}

  \begin{section}{Zero Forcing Sets and Maximum Nullity}\label{ZF}
  In this section we investigate a specific statistic on graphs, the zero forcing number.\\
  
  \begin{defn}Let $G=(V,E)$ be a graph and let $B\subseteq V$.  Color the vertices of $B$ black and color those not in $B$ white.  The \textbf{zero forcing rule} says that a black vertex $v$ can force a neighboring white vertex to become black if it is the unique white neighbor of $v$.  We say that $B$ is a \textbf{zero forcing set} for $G$ if repeated applications of the zero forcing rule eventually lead to all of $V$ becoming black.  It is clear that if $U$ is a zero forcing set for $G$, then so is any superset of $U$.  It is interesting to compute the minimum size of a zero forcing set, we call this minimum the \textbf{zero forcing number} for $G$ and denote it $Z(G)$.
  \end{defn}
  
  To see this more clearly we show an example with the following graph $G$.
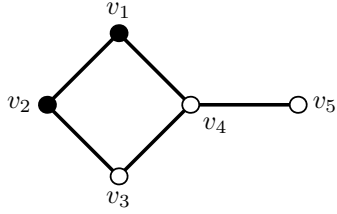
  \begin{center}
  \includegraphics{Graphics.3}
  \end{center}
  
  We see that $v_1$ can force $v_4$ because $v_1$ is black and $v_4$ is its only white neighbor.  $v_4$ cannot force anything yet because it still has two white neighbors.  However, $v_2$ can force $v_3$ and then $v_4$ can finally force $v_5$.  Now all of $V$ is black, and so $B=\{v_1, v_2\}$ was a zero forcing set.  This gives us that $Z(G)\leq 2$.  It is easy to see that there are no zero forcing sets of size 1, so, in fact, $Z(G)=2$.\\

This number is of interest for a variety of reasons; it has applications to quantum mechanics (see \cite{quantum}) but it appears more frequently in the mathematics literature for its relationship to the minimum rank problem (see eg. \cite{AIM08}, \cite{Barioli}).  Given a simple graph $G$, consider the set of symmetric matrices $\mathcal{S}$ such that the non-zero pattern of matrix in $\mathcal{S}$ agrees with the adjacency matrix of $G$ except, possibly, on the diagonal.  We allow the diagonal entries to be any element from the base field, including zero.  The problem is then to compute the minimum rank of a matrix in $\mathcal{S}$, denoted $mr^{\mathbb{F}}(G)$ where $\mathbb{F}$ is the base field.  Computing the maximum nullity for a matrix in the class, denoted $M^{\mathbb{F}}(G)$, is equivalent to computing the minimum rank since the rank-nullity theorem gives us $mr(G)+M(G)=|V|$.  It was proved in an AIM workshop (see \cite{AIM08}) that $Z(G)$ is an upper bound for the maximum nullity of $G$ over any field.\\
  
  \begin{thm}\label{ZgeqM}Let $G$ be a graph.  We then have $M^\mathbb{F}(G)\leq Z(G)$ for any field $\mathbb{F}$.
  \end{thm}

  Computing $M(G)$ exactly is equivalent to computing $mr(G)$, which we know is very hard in general.  Others have focused on using this bound to aid in computation of $M(G)$ for small graphs (see \cite{DeLoss}), but in this paper we focus on computing $Z(G)$ for bipartite circulants more generally.  We prove some upper and lower bounds on $Z(G)$ and discuss when they are actually attained.  We will start with the simplest lower bound, which depends only on the regularity of $G$.\\
  
  \begin{lem}\label{geqk}
  Let $G$ be a graph.
  \begin{enumerate}
    \item If $G$ is $k$-regular, then $Z(G)\geq k$.
    \item If $G$ is bipartite and $k$-regular, then $Z(G)\geq 2(k-1)$.
  \end{enumerate}
  \end{lem}

  \begin{proof}
  To see the first inequality, note that to get a black vertex to initiate a force it must have exactly one white neighbor.  Therefore, the first force requires that both the original vertex and $k-1$ of its neighbors are black for a total of at least $k$ black vertices.
  
  Analogously, if the graph is bipartite then the same analysis applies separately to each of the sets of the bipartition of the vertex set.  We need at least $2(k-1)$ black vertices, $(k-1)$ in each set of the bipartition, to get a force on each side.  Unlike in the non-bipartite case, if the forcing vertices are adjacent then we already counted them so we don't need to add them again, giving the result.
  \end{proof}
  
  A natural question to ask next is whether or not the lower bounds in the previous lemma can be achieved.  The first result is trivially achieved for the complete graph on $k+1$ vertices, $K_{k+1}$, and second by the complete bipartite graph with $k$ vertices in each partition, $K_{k,k}$, but there are certainly other examples in both cases as well.  We consider which bipartite circulants achieve the lower bound in the second case of the proposition, and we give a more general result which will give us some cases of equality as a corollary.\\
  
  \begin{thm}\label{leq2ik}
  Let $G$ be a connected bipartite circulant on $2n$ vertices.  If $G\simeq I+P^{i_2}+\cdots+P^{i_k}$ with  $0<i_2<\ldots < i_{k} \leq n-1$, then $Z(G)\leq 2i_k$.
  \end{thm}
  \begin{proof}
  We will exhibit a zero forcing set of $G$ with size $2i_k$.  The graph is the following (some unimportant edges not shown) with the first $i_k$ vertices on each side, corresponding to the rows and columns labeled 0 through $i_k-1$, initially colored black:
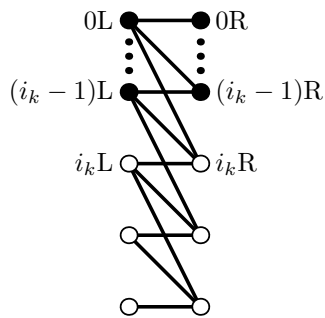
\begin{center}
\includegraphics{Graphics.4}
\end{center}
Now we can see that the first vertex on the left, 0L, has all black neighbors except for one because $i_k$ is the largest power of $P$.  It will force the first white vertex on the right, $i_k$R, which then has only one white neighbor and so it forces $i_k$L.  1L can now force $(i_k+1)$R which can force $(i_k+1)$L.  This will repeat itself moving down one row each time until every row has only black vertices.  This shows that we have a zero forcing set with $2i_k$ vertices.
  \end{proof}
  We have an immediate corollary.\\
  
  \begin{cor}\label{2(k-1)}
  Let $G$ be a connected bipartite circulant on $2n$ vertices.  If $G\simeq I+P+P^2+\cdots+P^{k-1}$ then $Z(G)= 2(k-1)$.
  \end{cor}
  \begin{proof}
  Since $G$ is bipartite and $k$-regular, we know from Lemma \ref{geqk} that $Z(G)\geq 2(k-1)$. Theorem \ref{leq2ik} gives us that $Z(G)\leq 2(k-1)$.  Thus $Z(G)=2(k-1)$.
  \end{proof}

  We might ask if there are other bipartite circulants that also achieve the lower bound.  The next result clarifies exactly when equality occurs for bipartite circulants.  To write this concisely, we need some more notation; we let $P^{[\alpha]_j}$ denote the sum $P^{\alpha}+P^{\alpha+j}+P^{\alpha+2j}+\cdots+P^{\alpha-j}$ so that it consists of a number of summands equal to the order of $j$ in $\mathbb{Z}/n\mathbb{Z}$.  When it is clear from context, we will omit the $j$.  In addition, $P^{[[\gamma]_c]_d}$ will denote the sum $P^{[\gamma]_d}+P^{[\gamma+c]_d}+P^{[\gamma+2c]_d}+\cdots +P^{[\gamma-c]_d}$ which consists of a number of terms (each of which is itself a sum) equal to the order of $c$ in $\mathbb{Z}/\text{gcd}(n,d)\mathbb{Z}$.\\

  \begin{thm}\label{lbeq}
  If $G$ is a connected, $k$-regular bipartite circulant with $2n$ vertices and $Z(G)=2(k-1)$, then its biadjacency matrix is isomorphic to one of the following matrices:
  \begin{enumerate}
    \item $I+P+P^2+\cdots+P^{k-1}$,
    \item $J-P^{\alpha}-P^{\alpha+i}-\cdots -P^{\alpha+ri}$ for some $\alpha$, $r$, and $i$ with gcd$(i,n)>1$ and with the order of $i$ in $\mathbb{Z}/n\mathbb{Z}$ greater than $r+1$.
    \item $P^{[[\gamma_1]_c]_d}+\cdots + P^{[[\gamma_{c-1}]_c]_d}+P^{\beta}+P^{[\beta+\alpha]_d}+\cdots + P^{[\beta+l\alpha]_d}$ where $d$ divides $n$, $c$=gcd$(d, \alpha)$, $l\geq 1$, and one of the following:
     \begin{enumerate}
       \item $\alpha \in [[\gamma_j]]$ for some $1\leq j\leq c-1$, or
       \item $\beta\in[-m\alpha]$, $l>m$ $($so $P^{[0]_d}$ and $P^{[\alpha]_d}$ are terms in the matrix$)$, and $c\leq m+2k$.
     \end{enumerate}
  \end{enumerate}
  Furthermore, each of these matrices gives rise to a bipartite circulant achieving $Z(G)=2(k-1)$ for every valid choice of parameters.
  \end{thm}
  We first note that this is an extension of Corollary \ref{2(k-1)} since it includes case (1).  Then $G\simeq I+P+P^2+\cdots+P^{k-1}$ and we know this graph achieves the lower bound for $Z(G)$ from the corollary.  This result is philosophically an expanded converse to that statement; we will assume $Z(G)$ equals the lower bound and deduce its structure.\\

\begin{proof}  
  Let $G$ be a bipartite circulant that satisfies the hypotheses and let its biadjacency matrix be written $P^{i_1}+P^{i_2}+\cdots+P^{i_k}$.  As in Lemma \ref{geqk}, a zero forcing set for $G$ of size $2(k-1)$ consists of two adjacent black vertices with each of their other $k-1$ neighbors also colored black.  Without loss of generality, we can assume that the forcing black vertices are in position $0$L and $0$R by shifting the labeling by a constant as in Theorem \ref{BCiso}; we also assume $i_1$=0 because it makes the proof a little easier to read.\\
  
  Either the vertices of $G$ are all black after the first force -- trivially satisfying (1) -- or we have at least one more force on both sides.  The vertices that are currently black are $\{0,-i_2,-i_3,\hdots,-i_k\}$L and $\{0,i_2,i_3,\hdots,i_k\}$R.  Without loss of generality (we won't require the powers to be increasing), say $-i_2$L forces.  Its neighbors are $\{-i_2,0,i_3-i_2,i_4-i_2,\hdots,i_k-i_2\}$R and all but one must be black.  Furthermore, if more than one vertex could force, we pick the one with the largest gcd($i_2,n)$.  This means that for all but one value of $l$, the equation\[i_l-i_2=i_m\mbox{ or equivalently } i_m+i_2=i_l\] has a solution for $m$.  This requires that all the powers differ from each other by $i_2$.  The proof now diverges into two cases.\\
  
  If the order of $i_2$ in $\mathbb{Z}/n\mathbb{Z}$ is greater than $k$ we can add a multiple of $i_2$ to each power by Theorem \ref{BCiso} to get $G\simeq I+P^{i_2}+P^{2i_2}+\cdots+P^{(k-1)i_2}$.  Since we assumed $G$ was connected and $i_2$ divides every power of $P$ in the biadjacency matrix, Lemma \ref{conn} gives us that gcd($i_2,n)=1$.  This implies that $i_2$ is a unit in the group and so we can use Theorem \ref{BCiso} to multiply all the powers by $i_2^{-1}$.  This lets us write the biadjacency matrix of $G$ as $G\simeq I+P+P^2+\cdots+P^{k-1}$.\\
  
  If the order of $i_2$ is less than or equal to $k$, let gcd$(i_2, n)=d$.  We can decompose $\mathbb{Z}/n\mathbb{Z}$ as residue classes mod $d$ and we see from the previous paragraph that we can write $G$'s biadjacency matrix as a sum of some complete residue classes and exactly one partial residue class.  Using the convention above we can write \[G\simeq P^{[\alpha_1]}+\hdots+P^{[\alpha_m]}+P^{\beta}+P^{\beta+i_2}+\hdots+P^{\beta+li_2}\] where $\beta$ and each $\alpha_j$ are members of distinct residue classes mod $d$.  Now we see that $-i_2$L is adjacent to one white vertex, $(\beta-i_2)$R, and it can be forced to become black.  Now $(-2i_2)$L is black and adjacent to only one white vertex in position $(\beta-2i_2)$R so it can also be forced.  This repeats (on both sides) until we have $\{[-\alpha_1],[-\alpha_2],\hdots,[-\alpha_m],[-\beta]\}$L and $\{[\alpha_1],[\alpha_2],\hdots,[\alpha_m],[\beta]\}$R as our set of black vertices.  If this is all of $G$ then we had a zero forcing set of minimum size and $G$ is almost a complete graph; it's only missing consecutive edges from one congruence class.  For ease of notation, we let $i_2=i$ in the description of the matrix.  Therefore we can write $G\simeq J_n-P^{\alpha}-P^{\alpha+i}-\cdots-P^{\alpha+ri}$ where $\alpha$ is in the partial congruence class, gcd$(i,n)>1$ and the order of $i$ is greater than $r+1$.  This is case (2).\\
  
  We are left with the case where the initial set of forces didn't color all of $G$ black.  So we currently have $m+1$ completely black residue classes in each set of the bipartition of the vertices.  Now if a left vertex's label is congruent to $j$ mod $d$, it is adjacent to all vertices congruent to $j+\alpha_1, j+\alpha_2,\hdots,j+\alpha_m$ and $l+1$ vertices in congruence class $j+\beta$.  Since the congruence classes have more than one element in them, each neighboring full class must already be black for vertex $j$ to force anything else, the partial class it's adjacent to must be white, and $j$ can have only one neighbor in that class.  This forces $l$ to be 0.    We can therefore write \[G \simeq P^{[\alpha_1]}+\hdots+P^{[\alpha_m]}+P^{\beta}.\]
  
  The initial set of forces gives leaves us with the vertices $\{[\alpha_1], [\alpha_2], \hdots, [\alpha_m], [\beta]\} $R and $\{[-\alpha_1], [-\alpha_2], \hdots, [-\alpha_m], [-\beta]\}$L colored black.  This is not all the vertices of $G$, so we have at least one more force on each side.  Without loss of generality, assume that $-\alpha_1$L forces (relabel the class with the representative that forces if necessary).  Then, because the black vertices comprise entire residue classes, it must force along the $P^{\beta}$ edge.  This requires that $[\alpha_i-\alpha_1]$R is black for each $i$, and $(\beta-\alpha_1)$R is white.\\
  
  The same reasoning that we used for $i_2$ applies here to the classes; every class of edges mod $d$ except one differs from another by $\alpha_1$.  Let $c =$ gcd($d,\alpha_1)$.  Then the classes are of size $\frac{d}{c}$ and the $\alpha_i$ break up into one partial class mod $c$ and potentially several complete double classes (by double, we mean mod $c$ within mod $d$) as well.  Since we know that $(\beta-\alpha_1)$R is white, the $\beta$ edge must be part of the partial class mod $c$.  Combining this information into one statement and writing $\alpha_1$ as $\alpha$ for notation's sake, gives us that \[G\simeq P^{[[\gamma_1]_c]_d}+\cdots+ P^{[[\gamma_k]_c]_d}+P^{\beta}+P^{[\beta+\alpha]}+\cdots+P^{[\beta+l\alpha]}.\]
  
  We now show $l\geq 1$.  Assume not, then $\beta$ is a partial class by itself.  Since $\beta$ is the only edge of $G$ not included in a class, we could have written the double classes of edges as single classes mod $dc$ rather than as a sum of residue classes mod $d$ and then mod $c$.  This contradicts our assumption that we took the $i_2$ with the largest gcd originally, so this cannot occur and $l\geq 1$.\\
  
  Furthermore, we have at least two edges of $G$ that belong to each of $[[\gamma_i]]$ and $[[\beta]]$.  So even if the rest of $[[\beta]]$R gets forced, we will never be able to force a vertex in a new double class (not one of the $[[\gamma_i]]$ or $[[\beta]]$) since every other double class is currently entirely white.  Therefore, we must have $c-1$ complete double classes if we are going to force all of $G$, and this lets us further refine the structure of $G$:\[G\simeq P^{[[\gamma_1]_c]_d}+\cdots+ P^{[[\gamma_{c-1}]_c]_d}+P^{\beta}+P^{[\beta+\alpha]}+\cdots+P^{[\beta+l\alpha]}.\]
    
We must have $\alpha$ as an edge of $G$ to get $-\alpha$L to force (which we assumed), so we are in one of two cases: $\alpha$ belongs to a $[[\gamma_i]]$ or $\alpha$ belongs to $[[\beta]]$.  If $\alpha\in [[\gamma_i]]$ then every vertex on the right outside of $[[\beta]]$R is black and $[-\alpha]$L can force $[\beta-\alpha]$R (one vertex at a time).  Since $\alpha$ is part of a complete double class, $[-2\alpha]$L is also black and it can then force $[\beta-2\alpha]$R.  This continues until all of $[[\beta]]$R is forced; now all of the vertices are black (once we do the symmetric forces from R to L as well) and so this is a valid case, it's case (3a).\\
  
If $\alpha \in [[\beta]]$ then we can rewrite $G$ as \[G\simeq P^{[[\gamma_1]_c]_d}+\cdots+ P^{[[\gamma_{c-1}]_c]_d}+P^{\beta}+P^{[(-m+1)\alpha]}+\cdots+P^{[0]}+\cdots +P^{[k\alpha]}.\]  We assumed that $-\alpha$L forced, so we must have $k\geq 1$.  Now $[-\alpha]$L forces $[\beta-\alpha]$R (one vertex at a time).  Next $[-2\alpha]$L can force $[\beta-2\alpha]$R as long as $k\geq 2$.  This continues until $[-k\alpha]$L forces $[\beta-k\alpha]$R.  The set of currently black vertices is $\{[[\gamma_1]],\ldots,[[\gamma_{c-1}]],[\beta-k\alpha],\ldots,[k\alpha]\}$ and so if $c\leq 2k+m$ we have forced all of $[[\beta]]$R and this is a zero forcing set for $G$.  This is case (3b).\\

If not, we have that both $(-k-1)\alpha$L and $(\beta-k-1)\alpha$R are white.  Since the only edge that can generate forces is $\beta$, these two vertices will always be white.  $G$ does not satisfy $Z(G)=2(k-1)$ and so $G$ must have been one of the above cases.
  \end{proof}
  
  Note that the most complicated case of the previous theorem can only occur when there are at least three factors in $n$'s prime factorization.  If $n$ is prime the result condenses considerably, leading to the following corollary.
  
  \begin{cor}
  If $G$ is a connected, $k$-regular bipartite circulant on $2n$ vertices, $n$ is prime, and $Z(G)=2(k-1)$, then $G\simeq I+P+P^{2}+\cdots+P^{k-1}$.
  \end{cor}
  \begin{proof}
  In Theorem \ref{lbeq} cases (2) and (3) can only happen when a power of $P$ has a non-trivial gcd with $n$.  Since $n$ is prime, this cannot happen.  We must therefore be in case (1).
  \end{proof}
  
  The regularity of the graph is also smallest in case (1) of Theorem \ref{lbeq}, so it gives us the best bound on $mr(G)$.  We state the corresponding result as a corollary.\\
  
  \begin{cor}
  If $G$ is a bipartite circulant on $2n$ vertices and $G\simeq I+P+P^2+\cdots+P^{k-1}$, then $mr(G)\geq 2n-2k+2$.
  \end{cor}
 \begin{proof}
 From Theorem \ref{lbeq} we see that case (1) applies and $Z(G)=2(k-1)$.  From Theorem \ref{ZgeqM} we see that $M(G)\leq 2(k-1)$ so $mr(G)\geq 2n - 2(k-1)$ by the rank-nullity theorem.
 \end{proof}
 
  In the next section we show that both the upper and lower bound in this section can be improved if we restrict the regularity of the graph.
  
  \end{section}
  
  \begin{section}{Zero Forcing in 3-Regular Bipartite Circulants}\label{cubic}
  Exact computation of $Z(G)$ becomes more difficult for bipartite circulants as the degree of regularity increases.  Because of this, we can extend the results of Section \ref{ZF} by restricting ourselves to the first interesting case -- $G$ is a cubic bipartite circulant.  We will use this restriction to extend the characterization of equality in Theorem \ref{lbeq} to a much stronger result and then show a more explicit structural result in some cases and use it to obtain another set of bounds for $Z(G)$.\\
  
We can extend Theorem \ref{lbeq} in the case where $G$ is 3-regular and bipartite.  Note that we have lost the bipartite circulant assumption and are in fact characterizing the unique cubic bipartite graph achieving $Z(G)=4$.
  
  \begin{thm}
  If $G$ is a connected, bipartite, cubic graph on $2n$ vertices with $Z(G)=4$, then $G\simeq I+P+P^2$.
  \end{thm}
  \begin{proof}
  We need many pictures in this proof and we adopt the convention that each figure has all the currently black vertices present.  More vertices could be forced as the argument in each case progresses, but we attempt to infer the structure of $G$ from its zero forcing set.\\
  
  From Lemma \ref{geqk} it is clear that we need at least two black vertices in each bipartition in order to obtain a force on each side.  Without loss of generality, after one force on each side we can assume we have situation drawn below.
  
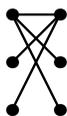
\begin{center}
  \includegraphics{IPP2.0}
  \end{center}
  We must examine the bottom four vertices.  We will look at how many connecting edges there are and do each case separately.  Consider all six vertex subgraphs corresponding to a zero forcing set of size 4 as pictured above.  We assume that we have a zero forcing set that generates a six vertex subgraph with the most edges, we denote the number of edges beyond the required six by $e$.  The proof then breaks up into cases based on the value of $e$.\\

If $e$ is zero, no other forces will be possible so $Z(G)\neq 4$.  Therefore $e \geq 1$.  If $e$ is one then the graph looks like this:
  
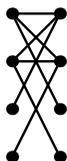
\begin{center}
  \includegraphics{IPP2.1}
  \end{center}
  Now if we have any edges between the new bottom four vertices we could find a different six vertex subgraph with the same zero forcing set that has more than one edge.  If we take the vertices of the first and second rows and both vertices adjacent to the extra edge then we get a six vertex subgraph of $G$ of the appropriate type giving $e\geq 2$.  If we don't have any connections between the bottom four vertices then we can't get any more forces since every black vertex has either zero or two white neighbors.  This implies $Z(G)\neq 4$.  Therefore we must have $e\geq 2$.\\

  If $e$ is two then we have to consider two cases, either one of the four vertices has 3 black neighbors or each of them has exactly two black neighbors.  If a vertex had three black neighbors, without loss of generality, the graph would look like this:
  
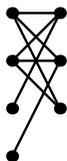
\begin{center}
  \includegraphics{IPP2.21}
  \end{center}
  Now either the third vertex on the right connects to the fourth on the left or it doesn't.  If it doesn't, we can't get any more forces on the right, since no two of the black vertices on the right connect to the same vertex with a remaining connection on the left.  If instead 3R is adjacent to 4L, then 2R, 1R, 3R and 1L, 2L, 4L form a six vertex subgraph of $G$ corresponding to the initial zero forcing set with three edges.  This is a contradiction.\\

  We are left with the $e=2$ case where none of the bottom four vertices has more than two neighbors within the subgraph.  Without loss of generality, we have:
  
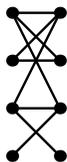
\begin{center}
  \includegraphics{IPP2.22}
  \end{center}
  Consider the 2L, 2R, 4L, 4R vertices (the ones that need higher degree than they currently have).  Then we have four vertices, each of which is black, with one edge between them.  Since we know we need at least two edges to get a graph with $Z(G)=4$, we're back in the original scenario but with two additional vertices.  We can use this case again to add an edge between the degree one vertices, or we could use an $e \geq 3$ case to add some edges and continue building $G$.  We can repeat this step connecting our bottom two vertices to two new vertices -- extending the ladder pattern -- but we eventually need to do something different with our four vertices in order to finish the graph.  We'll revisit this after considering the other cases.\\

If $e$ is three then without loss of generality we have:
  
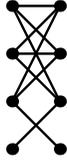
\begin{center}
  \includegraphics{IPP2.31}
  \end{center}
  
  Since the only vertices with degree less than three are the two in the bottom row, either they are adjacent or they are not.  If not, neither of them can have two black neighbors and so we can't get any more forces.  This eliminates this case.  If they are adjacent, then the graph is the above but with another row of vertices as follows:
  
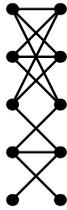
\begin{center}
  \includegraphics{IPP2.32}
  \end{center}
  Now the exact same analysis applies and we must have an edge between the vertices in the new bottom row.  However, they are required to have degree three, so we must continue adding more vertices to the graph.  Unfortunately, there is no way to ever finish this process with a 3-regular graph without requiring another vertex in the zero forcing set.  This implies $Z(G)\neq 4$, so we can eliminate this case as well.\\

  If $e$ is four, then every edge is present and $G$ is isomorphic to $K_{3,3}$ which has biadjacency matrix $I+P+P^2$ as desired.\\

  Therefore, the only way to extend the number of vertices in the graph is to repeat the second case for two connections, adding one vertex to each side for each repeat and increasing the size of the ladder.  In addition, the only way to end the process is to add all four of the potential connections ($K_{3,3}$ if the extending process isn't used).  This gives rise to graphs of the form:
  
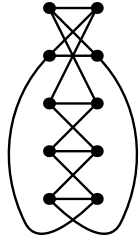
\begin{center}
  \includegraphics{IPP2.4}
  \end{center}
  The biadjacency matrix of this graph can be written as $I+P^{-1}+P^1$ (to see this, slide the second row above the top row).  Then $P(P^{-1}+I+P^1)=I+P^1+P^2$, giving the result.
  \end{proof}
  
  We might ask whether the corresponding result is true for 4-regular bipartite graphs $G$ with $Z(G)=6$.  Unfortunately there are many 4-regular bipartite graphs with $Z(G)=6$, not all of them are isomorphic. With $n=8$, $I+P+P^2+P^3$ has rank 5 and we know it's the only 4-regular bipartite circulant in $n=8$ with $Z(G)=6$ from Theorem \ref{lbeq}.  However, the matrix
  \begin{center}
  \[ A= \left( \begin{array}{cccccccc}
1 & 1 & 0 & 1 & 1 & 0 & 0 & 0\\
1 & 1 & 1 & 1 & 0 & 0 & 0 & 0\\
1 & 1 & 1 & 0 & 0 & 1 & 0 & 0\\
1 & 0 & 1 & 0 & 0 & 1 & 1 & 0\\
0 & 1 & 0 & 0 & 0 & 1 & 1 & 1\\
0 & 0 & 1 & 1 & 1 & 0 & 0 & 1\\
0 & 0 & 0 & 1 & 1 & 0 & 1 & 1\\
0 & 0 & 0 & 0 & 1 & 1 & 1 & 1\end{array} \right)\] 
  \end{center}
  has four ones in each row and column, rank 6, and $G\simeq A$ has $Z(G)=6$.  Therefore, $G$ cannot be expressed as a bipartite circulant since the rank is not changed by row and column permutations.\\
  
  We now investigate situations where we can get upper bounds that are smaller than those given in Theorem \ref{leq2ik}.  When one or more of the powers of $P$ in the biadjacency matrix have non-trivial gcds with $n$, it is difficult to use Theorem \ref{BCiso} to reduce the size of the powers because multiplying by a unit will never reduce that power below the gcd.  As we saw in the proof of Theorem \ref{lbeq}, if $G$ is more than three regular with non-trivial gcds with $n$ and the powers of $P$ it is difficult to analyze the structure of $G$.  However, if we restrict ourselves to cubic bipartite circulants we can continue our discussion.  This will give us a structural result that we can use to bound $Z(G)$ for such graphs.\\
  
      \begin{lem}\label{gcd}
    Let $G$ be a connected, cubic bipartite circulant on $2n$ vertices with $G\simeq I+P^i+P^j$ and gcd$(j,n)=d>1$.  Then $G$ decomposes into $d$ cycles of length $\frac{2n}{d}$, each of which is connected to neighboring cycles by a matching on half of the vertices of each cycle.
    \end{lem}
    \begin{proof}
    In Section \ref{BC} we saw that the subgraph of $G$ with the $P^i$ edges removed is a set of $d$ disjoint cycles each of length $\frac{2n}{d}$ with each cycle corresponding to a different number mod $d$.  Consider the cycle with 0L.  The $P^i$ edges goes from 0L to $iR$, $j$L to $(j+i)$R, $2j$L to $(2j+i)$R, and so on.  Since $d$ divides $j$, every multiple of $j$ plus $i$ is in the same congruence class, and so is in the same cycle.  Since $G$ is connected, the cycles are different and we have a matching from the left side of the congruent-to-zero cycle to the right side of the congruent-to-$i$ cycle.  The left side of the congruent-to-$i$ cycle is similarly matched with the right side of the congruent-to-$2i$ cycle and so on.  Since gcd$(i,j,n)=1$ by Lemma \ref{conn}, we see that every one of the cycles is eventually connected by a multiple of $i$.  This gives our result.
    \end{proof}

  Now, if $G$ is a cubic bipartite circulant and one of the powers in its matrix representation is not relatively prime to $n$, we can use Lemma \ref{gcd} to get both an upper and a lower bound for $Z(G)$.\\
  
  \begin{thm}\label{Zgcd}
  Let $G$ be a connected, cubic bipartite circulant on $2n$ vertices with $G\simeq I+P^i+P^j$ and gcd$(j,n)=d>1$.  Then $d < Z(G) \leq d+ 2(\frac{n}{d})-2$.
  \end{thm}
  \begin{proof}
  From Lemma \ref{gcd} we know that $G$ decomposes into $d$ cycles $C_0,\ldots, C_{d-1}$ of length $\frac{2n}{d}$ which are connected to each other by edges from $P^i$.  We label each of the cycles by what each vertex in them is congruent to mod $d$.  Furthermore, the L vertices in $C_m$ are matched with the R vertices from $C_{m+i}$ (reduce mod $d$ if necessary) by $P^i$ so that each cycle $C_m$ is adjacent to exactly two other cycles, $C_{m+i}$ and $C_{m-i}$, via edges from $P^i$.\\
  
  To obtain the lower bound we consider the vertex sets $V_m$ consisting of the L-vertices in $C_m$ and the R-vertices in $C_{m+i}$.  There are $d$ such sets, one for each choice of $m\in \mathbb{Z}/d\mathbb{Z}$.  Now every $v\in V_m$ is adjacent to one vertex in $V_m$ along its edge from $P^i$ and two vertices outside the set.  The two neighbors of $v$ outside $V_m$ are from either $C_m$ or $C_{m+i}$ depending on whether $v$ was an L-vertex or an R-vertex, but in either case they are adjacent to $v$ and another vertex in $V_m$ by following edges in the cycle.  Now if a coloring of $V$ had all of the vertices in $V_m$ white, none of the vertices could ever be forced because they could never be the sole white neighbor of an outside vertex.  This requires any zero forcing set for $G$ to contain at least one vertex from each of these vertex sets, so $Z(G)\geq d$.  If we had exactly one black vertex per set we couldn't have a black vertex with two black neighbors, so we'd need more black vertices to get any forces.  This proves the lower bound.\\
  
  To obtain the upper bound we exhibit a zero forcing set with the appropriate number of vertices.  First we note that if cycle $C_m$ is black, it can force all of the R vertices of $C_{m+i}$ (the rest of $V_m$) and all of the L-vertices of $C_{m-i}$ (the rest of $V_{m-i}$).  If even one of the L-vertices of $C_{m+i}$ is black, the rest of the cycle will be forced from vertices in $V_m$.  This would allow us to repeat the argument and force the rest of $V_{m+i}$.  Therefore if we start with all of $C_0$ black ($2(\frac{n}{d}$) vertices) and one L-vertex from each of the $V_m$ other than 0 and -$i$ ($d-2$ vertices) we can force the rest of $V_0$.  Along with the one black vertex on the other side of $C_i$ this forces the rest of $C_i$ which can force the rest of $V_i$.  We continue in this fashion until we get to $C_{-i}$ which doesn't have an L-vertex already black.  However the half in $V_{-i}$ is already black from the initial forces, and so we can force every vertex.  Thus we have a zero forcing set of the required size.  This proves the upper bound.
  \end{proof}
  
  In fact the upper bound is actually achieved in several examples.  For example, if $G$ is a connected bipartite circulant on 2$n$ vertices with $G\simeq I+P+P^{\frac{n}{2}}$, it is true that $Z(G)=\frac{n}{2}+2$ for $n\geq 4$.\\
  
  \end{section}
  \begin{section}{Conjectures and Future Research}\label{conj}
  If we want to use the bounds we found in Section \ref{ZF} to get the smallest upper bound on $Z(G)$ for a bipartite circulant, we need to know which sets of powers of $P$ generate a biadjacency matrix for $G$.  We saw in Theorem \ref{BCiso} that we can often find several ways of representing a given bipartite circulant graph with circulant biadjacency matrices.  This leads us to ask if there are any other potential representations of the biadjacency matrix in circulant form, and the following conjecture attempts to answer that question.\\
  
  \begin{conj}
  Let $G$ be a connected bipartite circulant on $2n$ vertices with both $G\simeq P^{i_1}+\cdots + P^{i_k}$ and $G\simeq P^{j_1}+\cdots + P^{j_k}$.  Then there exists a $u\in \mathbb{Z}/n\mathbb{Z}$ with gcd$(u,n)=1$ and $z\in \mathbb{Z}/n\mathbb{Z}$ such that $\{ui_1+z,\ldots,ui_k+z\}=\{j_1,\ldots,j_k\}$.
  \end{conj}
  
  In computations for small values of $n$ this appears to be valid for cubic bipartite circulants, but it is difficult to compute in general.\\
  
  Computations reveal that the upper bounds from Theorem \ref{leq2ik} and Theorem \ref{Zgcd} are attained for cubic bipartite circulants with relatively small values of $n$.  However, to be most useful, they both require that we compute all potential circulant bidajacency matrices realized by $G$ using Theorem \ref{BCiso}.  This brings us to our next question:\\
  
  \begin{q}
  Let $G$ be a connected cubic bipartite circulant on $2n$ vertices.  Under what conditions on $n$ and the powers of $P$ in the matrix is the zero forcing set given by Theorem $\ref{leq2ik}$ minimal? Similarly, when is the bound from Theorem $\ref{Zgcd}$ minimal?  Note that we may need to consider all circulant biadjacency matrix representations of $G$ because the bounds depend on the specific circulant matrix, but $Z(G)$ depends only on $G$.
  \end{q}
  
  One of the natural questions to ask about any statistic that bounds the minimum rank of a graph is when the actual minimum rank is equal to the statistic.  For an arbitrary graph this is often difficult, but it may be easier for bipartite circulants because of the high degree of predictability in their adjacency matrices.\\
  
  \begin{q}\label{MeqZ}
  When is $Z(G)=M(G)$ for bipartite circulant graphs?
  \end{q}
  
  If Question \ref{MeqZ} proves intractable, can we at least bound $M(G)$ above for a given $n$ and degree of regularity?  In other words,\\
  
  \begin{q}
  Is there a simple function that depends only on $n$ and $k$ such that $Z(G)\leq f(n,k)$ for all connected $k$-regular bipartite circulants on $2n$ vertices?
  \end{q}
  
  This would give us a lower bound on the minimum rank over the whole class of bipartite circulants, which is difficult in general.

 \end{section}
\bibliographystyle{plain}

  \end{document}